Article: CJB/2010/5

# Omega Circles

**Christopher J Bradley**

**Abstract:** Circles through the Brocard points (Omega circles) share nearly all the properties of circles through the orthocentre including the fact that key triangles inscribed in them are indirectly similar to triangles inscribed in the circumcircle. Properties of Omega circles are described in this article, thereby concluding our work on indirectly similar triangles. It is also shown that the three points where the medians intersect the orthocentroidal circle are such that circles through these points carry triangles directly similar to triangles inscribed in the circumcircle.

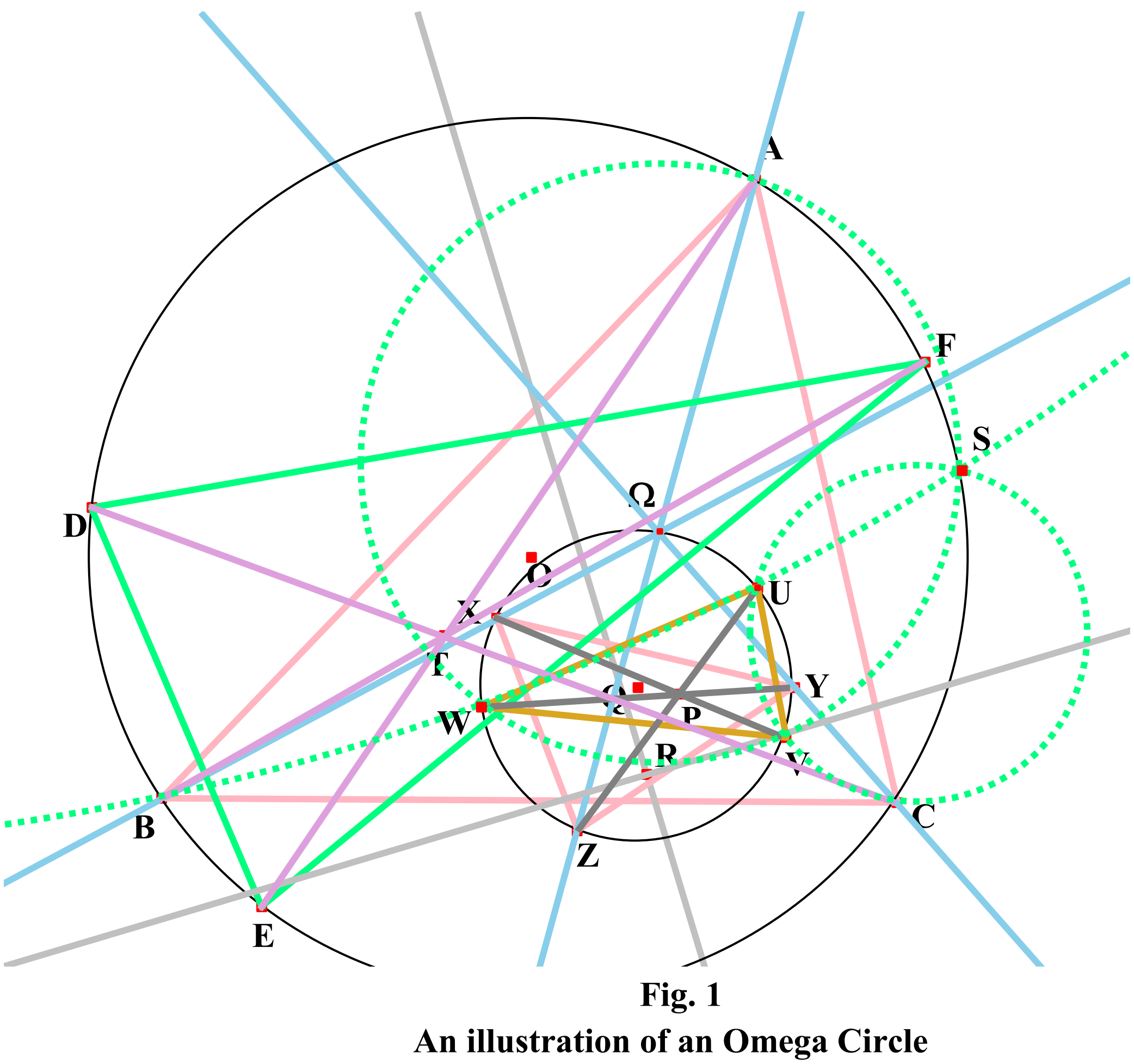

**Fig. 1**
**An illustration of an Omega Circle**

## 1. Introduction

Every circle passing through H, the orthocentre of a triangle ABC is called a Hagge circle, named after the first person [1] to detail their most important properties. These are (i) If AH, BH, CH meet a Hagge circle Σ at points X, Y, Z respectively, then triangle XYZ is indirectly similar to triangle ABC, (ii) If circles BHC, CHA, AHB meet Σ at U, V, W respectively, then UX, VY, WZ meet at a point P and (iii) P is also the centre of inverse similarity, meaning there is an axis through P such that if ABC is reflected in this axis followed by an appropriate enlargement (reduction) centre P, it is mapped on to triangle XYZ. It may be added that if triangle UVW is mapped back (by the inverse transformation) from Σ to the circumcircle to form a triangle DEF then DEF is inversely similar to UVW and because of property (ii) APD, BPE, CPF are straight lines.

If you take a general point J rather than H and repeat the construction then the triangle XYZ is no longer similar to ABC. That might have been the end of the story, but it now appears that circles through the Brocard points Ω and Ω' have properties nearly the same as Hagge circles, but with subtle differences. To be particular we restrict our discussion for the time being to Omega circles, those passing through the Brocard point Ω, which if AB > AC is further from BC than the Brocard point Ω' and which has areal co-ordinates $(1/b^2, 1/c^2, 1/a^2)$. The properties of an Omega circle Γ are (i) If AΩ, BΩ, CΩ meet Γ again at Z, X, Y respectively, then triangle XYZ is indirectly similar to triangle ABC, (ii) If circles BΩC, CΩA, AΩB meet Γ at U, V, W then UZ, VX, WY meet at a point P. However, the centre of inverse similarity R is in general distinct from P, so that when the points D, E, F are obtained by finding the inverse images of U, V, W respectively, the lines AE, BF, CD meet at a point T which is the inverse image of P under the indirect similarity centre R. The exception is the seven-point circle, when R and P coincide. These properties are all exhibited in Fig. 1, which is drawn using *CABRI II* plus. There is good reason to believe that H, Ω, Ω' are the only three points that exhibit the indirect similarities, (see Article: CJB/2010/4 of this series) but there are three other points which provide direct similarities rather than indirect similarities between ABC and XYZ. There are numerous other properties that emerge from both Hagge circles and Omega circles such as the fact that pairs of triangles are orthologic and that the perspectives also create axes of importance. In this article, we do not attempt to catalogue all these properties. The interested reader is referred to the work of Speckman [2].

## 2. The equation of an Omega circle

We use areal co-ordinates with A, B, C the triangle of reference. The point Ω has co-ordinates $(1/b^2, 1/c^2, 1/a^2)$. From the equations of AΩ, BΩ, CΩ we may write down possible co-ordinates of Z, X, Y which are $Z(1, l/c^2, l/a^2)$, $X(m/b^2, 1, m/a^2)$, $Y(n/b^2, n/c^2, 1)$. However, l, m, n are related to one another because the circle XYZ passes through Ω. The condition for this is

$l = -\ a^2b^2c^2(a^2(b^2c^2 + n(c^2 - m)) + c^2m(b^2 + n))$ all divided by
$a^4b^2(b^2c^2 + n(c^2 - m)) + a^2c^2b^4(m - n) - b^2m(c^2 + n) - c^2mn - b^4c^2mn.$ (2.1)

What Equation (2.1) implies is that there are a doubly infinite number of Omega circles given by varying the values of m and n. By eliminating l in favour of m and n means that symmetry is lost. Despite that further working becomes easier to check using a programme such as *DERIVE*, which is the one we use. The equation of the Omega circle, in terms of m and n, is

$$b^4x^2(a^2(b^2c^2 + n(c^2 - m)) + c^2m(b^2 + n)) + b^2x(y(a^2(b^4c^2 - b^2(c^4 + mn) - c^4n) + c^2m(b^4 - b^2c^2 - c^2n)) + z(a^4n(b^2 + m) - a^2(b^4n + b^2m(c^2 + n) + c^2mn) - b^4mn)) + a^2(c^2my^2(b^2(c^2 + n) + c^2n) - c^2yz(a^2(b^4 + b^2(m + n) + mn) + m(b^4 - b^2c^2 - c^2n)) + b^2nz^2(a^2(b^2 + m) + b^2m)) = 0. \qquad (2.2)$$

## 3. The indirect similarity

The fact that triangles ABC and XYZ are indirectly similar follows by angle chasing as a result of the facts that angle BΩC = $180^0 - C$, angle CΩA = $180^0 - A$, and angle AΩB = $180^o - B$. (If Ω is replaced by H then C, A, B are replaced by A, B, C in these equations, involving a *cyclic* rearrangement of A, B, C.)

## 4. The circles BΩC, CΩA, AΩB

The circle BΩC has equation

$$b^2x^2 - xy(c^2 - b^2) + a^2yz = 0. \qquad (4.1)$$

The equations of circles CΩA, AΩB may be written down by cyclic change of x, y, z and a, b, c.

## 5. The points U, V, W and P

If a triangle ABC is given and a point J and any circle Σ through J, it is a general result that if AJ, BJ, CJ meet Σ again in points Z, X, Y and circles BJC, CJA, AJB meet Σ again in points U, V, W then XV, YW, ZU are concurrent at a point P, so we do not repeat the proof, but simply give the co-ordinates of the points in the present context. (The result can then be checked in this case.) Note that the perspective is quite independent of the similarity in this case between ABC and XYZ.

U has co-ordinates (x, y, z) where

$$x = a^2m(b^2c^2 + n(b^2 + a^2)),$$
$$y = b^2n(a^2b^2 + m(a^2 + c^2)), \qquad (5.1)$$
$$z = (b^2m(a^2(b^2n + c^2(m - n) + 2mn) + mn(b^2 - c^2))(b^2(c^2 + n) + c^2n))/(n(a^2(b^2 + m) + b^2m)).$$

V has co-ordinates (x, y, z) where

$x = -(n(a^2(b^2+m)+b^2m)(a^4(b^2c^2+n(c^2-m))-a^2c^2(b^2(c^2-m+n)+n(c^2-m))-c^2m(b^2(c^2+n)+c^2n)))/(b^2(a^2(b^2c^2+n(c^2-m))+c^2m(b^2+n)))$,

$$y = n(a^2(b^2+m)+b^2m), \quad (5.2)$$

$$z = a^2(b^2c^2+n(c^2-m))+c^2m(b^2+n).$$

W has co-ordinates (x, y, z) where

$$x = a^2m(b^2(c^2+n)+c^2n),$$

$y = (b^2(a^2c^2(b^4+b^2(m+n)+mn)-b^4mn)(a^2(b^2c^2+n(c^2-m))+c^2m(b^2+n)))/(c^2m(b^2(c^2+n)+c^2n))$, (5.3)

$$z = b^2(a^2(b^2c^2+n(c^2-m))+c^2m(b^2+n)).$$

P has co-ordinates (x, y, z) where

$$x = (c^2mn(a^2(b^2+m)+b^2m))(b^2(c^2+n)+c^2n),$$

$$y = (b^2n(a^2(b^2+m)+b^2m))(a^2(b^2c^2+n(c^2-m))+c^2m(b^2+n)), \quad (5.4)$$

$$z = (b^2(a^2(b^2c^2+n(c^2-m))+c^2m(b^2+n)))(m(b^2(c^2+n)+c^2n)).$$

It may be proved that VW^AΩ, WU^BΩ, UV^CΩ, P are collinear, but that we leave to the reader. As with all indirect similarities the midpoints of AX, BY, CZ, DU, EV, FW lie on a conic.

## 6. The points R, T, S

The points R, T, S are the additional points in the figure.

In view of the indirect similarity between ABC and XYZ there are bound to be a point R and two axes through R that are the centre and axes of inverse similarity. Nothing additional is proved by working out the co-ordinates of R, as it is known to exist. Its position is constructed as follows. The directions of the axes are parallel to the bisectors of any pair of corresponding sides such as BC and YZ. Take any such line and reflect triangle ABC in it. It will now be found that AX, BY, CZ are concurrent at a point that does not lie on the line. Move the line parallel to itself and when the point lies on the line that will be R, the centre of inverse similarity. This construction, of course, cannot be done with complete precision. But the construction can be made precise by noting that R must lie on certain lines. For example the point Ω lies on the circle XYZ. Its angular position relative to these points can be mirrored precisely on the circumcircle of ABC. Call that point $\Omega_o$. Call the same Brocard point of triangle XYZ the point $\Omega_1$. The passage from $\Omega_o$ to $\Omega_1$ is now the product of two indirect similarities and is thus an enlargement, so that $\Omega_o\Omega_1$ must pass through R. A second such point, leading to a second line through R, might be the second point of intersection of the circle XYZ with the seven-point circle. The two lines fix the position of R exactly.

Since circles BUC, CVA, AWB are concurrent at Ω, it follows that circle AVW, BWU, CUV are concurrent at a point S. As indicated by *CABRI II* plus, S lies on the circumcircle of ABC. This curious fact we have been unable to prove, as *DERIVE* found the equations too technically complicated to solve.

## 7. Omega prime circles

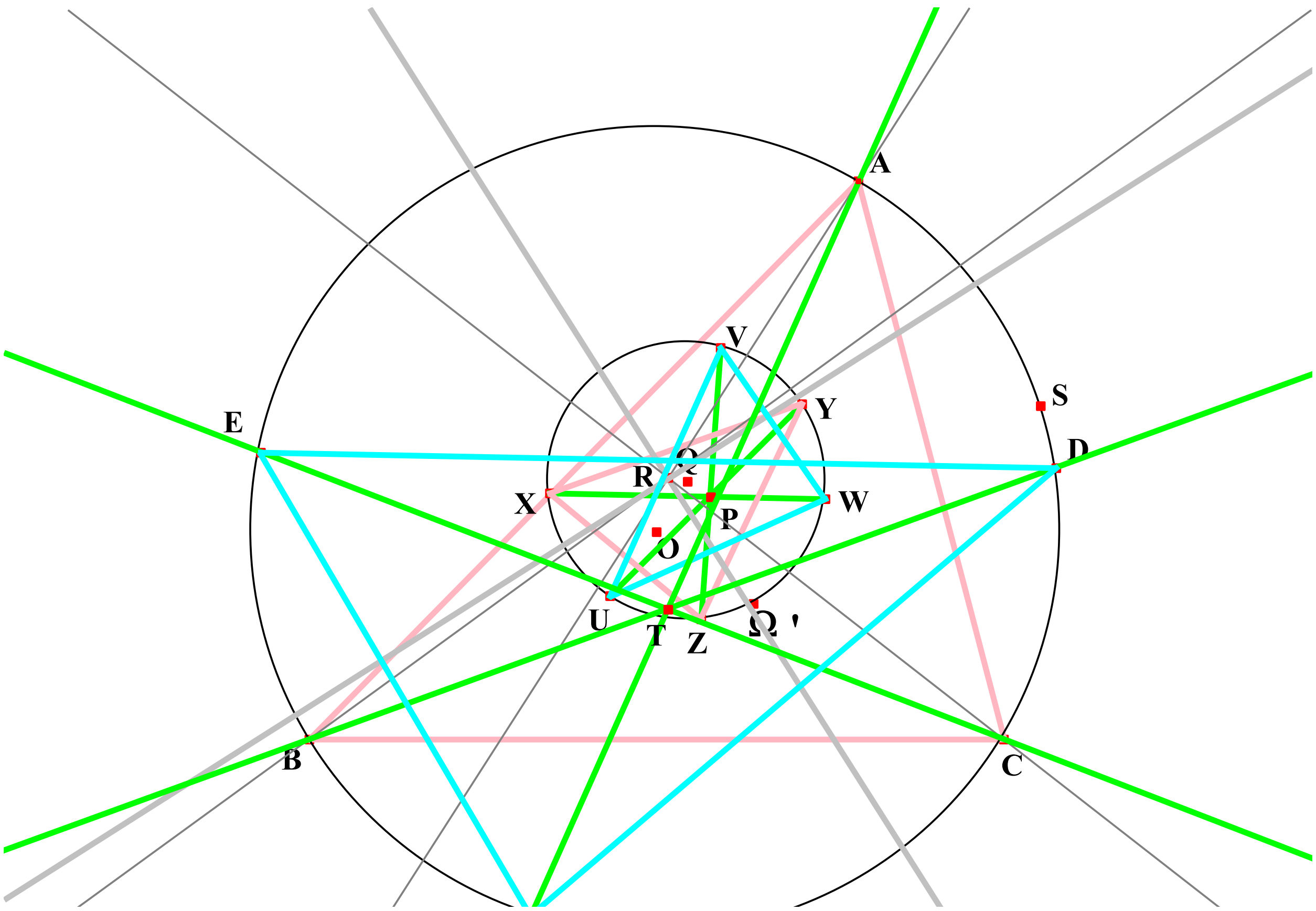

**Fig.2**
**An illustration of an Omega prime circle**

Circles through Ω' $(1/c^2, 1/a^2, 1/b^2)$ have very similar properties to those passing through Ω, the difference arising from the labelling of the points X, Y, Z. Now AΩ', BΩ', CΩ' meet the Omega prime circle at points Y, Z, X respectively. Then there is an indirect similarity between triangles ABC and XYZ. This similarity follows again by simple angle chasing bearing in mind that angle BΩ'C = $180^o$ – B, angle CΩ'A = $180^o$ – C and angle AΩ'B = $180^o$ – A. See Fig. 2.

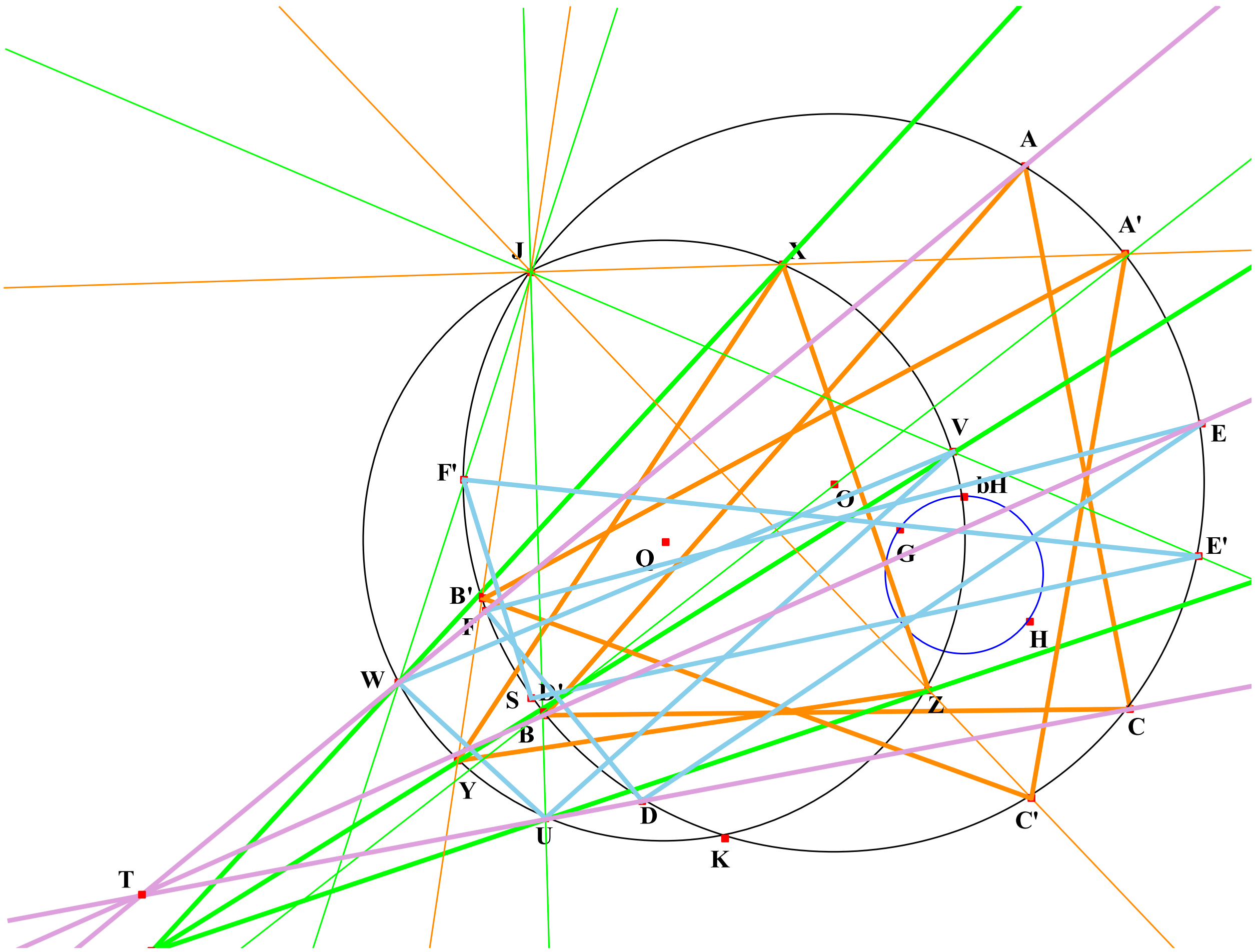

**Fig. 3**
**An illustration of a bH circle**

## 8. aH bH cH circles

aH, bH cH are the intersections of the lines AG, BG, CG with the orthocentroidal circle on GH as diameter, where G is the centroid of ABC and H is its orthocentre. We illustrate in Fig. 3 what happens when a circle Σ is drawn through bH. The lines AbH, BbH, CbH meet Σ again at Z, Y, X respectively. Since angle BbHC = $180^o$ – C, angle CbHA = $180^o$ – A, angle AbHB = $180^o$ – B it follows by simple angle chasing that triangle XYZ is directly similar to triangle ABC.

The intersections of Σ with the circumcircle are denoted by J and K. According to the theory of Wood [3] either J or K may be used as a perspector to project triangle XYZ into a triangle A'B'C', where A' is the intersection of JX with the circumcircle and B' and C' are similarly defined. Triangle A'B'C' is directly similar to XYZ and so directly congruent to ABC. It may be seen in the figure that they are related by a rotation about O. If Σ is not large enough to intersect the circumcircle it may be enlarged about its centre Q until it does. All that happens is that an extra step is required to illustrate the direct similarity. Once this has been done the rest follows as with Omega circles. Points U, V, W are defined as intersections with Σ of the circles BbHC, CbHA, AbHC respectively and UZ, VY, WX intersect at P. Points D, E, F are the inverse images of U, V, W and DC, EB, FA intersect at T. The point S is the intersection of circles AYZ, BZX, CXY and it also lies on the circumcircle.

Flat 4,
Terrill Court,
12-14 Apsley Road,
BRISTOL BS8 2SP.